\documentclass{amsart}

\usepackage{latexsym}
\usepackage{amsmath}
\usepackage{amssymb}
\DeclareMathOperator{\tor}{Tor}

\begin{document}


{\theoremstyle{plain}%
\newtheorem{theorem}{Theorem}[section]  \newtheorem{corollary}[theorem]{Corollary}
  \newtheorem{proposition}[theorem]{Proposition}
  \newtheorem{lemma}[theorem]{Lemma}
  \newtheorem{question}[theorem]{Question}
  \newtheorem{conjecture}[theorem]{Conjecture}
}
{\theoremstyle{remark}
\newtheorem{fact}{Fact}
\newtheorem{remark}[theorem]{Remark}
\newtheorem{construction}[theorem]{Construction}
}

{\theoremstyle{definition}
\newtheorem{definition}[theorem]{Definition}
\newtheorem{example}[theorem]{Example}
}

\newcommand{\m}[1]{\marginpar{\addtolength{\baselineskip}{-3pt}{\footnotesize \it #1}}}

\newcommand{\comp}[1]{#1^c}
\newcommand{\numbc}{\mbox{\#$\operatorname{comp}$}}
\newcommand{\reg}{\operatorname{reg}}

\newcommand{\compGS}{\comp{G_{S}}}
\newcommand{\N}{\mathbb{N}}
\newcommand{\iv}{\iota}


\begin{center}
{\bf ON THE LINEAR STRAND OF AN EDGE IDEAL}
\vspace{.25cm}

Mike Roth\\
Department of Mathematics and Statistics \\
Queen's University \\ 
Kingston, ON K7L 3N6  Canada\\
{\tt mikeroth@mast.queensu.ca}\\
\vspace{.25cm}

Adam Van Tuyl \\
Department of Mathematical Sciences \\
Lakehead University \\
Thunder Bay, ON P7B 5E1 Canada \\
\tt{avantuyl@sleet.lakeheadu.ca}
\end{center}

\title[ON THE LINEAR STRAND OF AN EDGE IDEAL]{}
\keywords{graphs, resolutions, Betti numbers, edge ideals}
\subjclass[2000]{13D40, 13D02, 05C90}

\begin{abstract}  Let $I(G)$ be the edge ideal
associated to a simple graph $G$.
We study the graded Betti numbers that appear in
the linear strand of the minimal free resolution
of $I(G)$.
\end{abstract}

\maketitle


\section*{Introduction}
Let $R = k[x_1,\ldots,x_n]$ with $k$ a field.
If  $G$ is a simple graph 
(no loops or multiple edges, but not necessarily connected) with
edge set $E_G$ and vertex set $V_G = \{x_1,\ldots,x_n\}$,
then associated to $G$ is 
the monomial ideal $I(G) := \langle\{ x_ix_j ~|~ \{x_i,x_j\} \in E_G\}
\rangle$ in $R$
called the {\it edge ideal of $G$}. 
We would like to
establish a dictionary between algebraic invariants
of  $I(G)$  and numerical data associated to the graph $G$. 
Such a dictionary might then allow us to prove graph theoretical
results algebraically, and vice versa.  
Among the many papers that have studied the properties of edge ideals,
we mention \cite{Fr,HVT,HHZ,HHZ1,J,JK,K,SVV,V,V1,V2,Z}.

In this paper we investigate the relation between the graded Betti numbers 
arising from the resolution of $I(G)$ and numerical information
associated to $G$.
Recall that if $I$ is a homogeneous ideal, then 
we can associate to $I$ a minimal
graded free resolution of the form
\[
0 \rightarrow \bigoplus_j R(-j)^{\beta_{l,j}(I)}
\rightarrow \bigoplus_j R(-j)^{\beta_{l-1,j}(I)}
\rightarrow \cdots
\rightarrow  \bigoplus_j R(-j)^{\beta_{0,j}(I)}
\rightarrow I \rightarrow 0
\]
where $l \leq n$ and $R(-j)$ is the $R$-module obtained by shifting
the degrees of $R$ by $j$.  The number $\beta_{i,j}(I)$, the $ij$th 
graded Betti number of $I$, equals the number of generators
of degree $j$ in the $i$th syzygy module. Since $I(G)$ is a square-free
monomial ideal, our principal tool to study the numbers 
$\beta_{i,j}(I(G))$ will be Hochster's formula (see Propositions \ref{prop: hochster}
and \ref{prop: edgeidealbetti}).

Some relations between the numbers $\beta_{i,j}(I(G))$ and the properties
of $G$ have already been established.
For example, from the construction of $I(G)$ we see that
$\beta_{0,2}(I(G)) = \#E_G$, and $\beta_{0,j}(I(G)) = 0$ for $j \neq 2$.
Formulas for the Betti numbers of the first syzygy module 
are given in \cite{EV}.   In \cite{K} it is shown that 
$\beta_{i,2(i+1)}(I(G))$ with $i \geq 1$
is the number of induced subgraphs of $G$ consisting of exactly $i+1$
disjoint edges.  Zheng \cite{Z} showed that when $G$ is a tree,
the regularity of $I(G)$, a measure of the ``size'' of the resolution,
is also  related to the number of disjoint edges of $G$.
Many other such relations 
between $\beta_{i,j}(I(G))$ and $G$
can be found in the thesis of Jacques \cite{J}.

In this note we shall focus on the graded Betti numbers of
$I(G)$ that appear in the linear strand.
When $I$ is generated by
elements of degree at least $d$, then it can be shown that
$\beta_{i,j}(I) = 0$ for all $j < i+d$.  Thus, 
the first graded Betti number of interest is $\beta_{i,i+d}(I)$ for 
each $1 \leq i \leq l$.
The numbers $\beta_{i,i+d}(I)$ for $i \geq 0$ describe 
the {\it linear strand} of the minimal resolution of $I$, i.e.,
they count the number of 
linear syzygies appearing in the resolution.  It is then
of interest to understand what the linear strand can tell
us about $I$ (cf. \cite{EGHP,EG,RW}).  
Because $I(G)$ is generated by monomials of degree $2$,
we will be interested in the numbers $\beta_{i,i+2}(I(G))$
for $i \geq 0$.

The contributions of this note are as follows.  
If $G$ has no induced 4-cycles, then we provide an exact formula
for $\beta_{i,i+2}(I(G))$ in terms of data associated to $G$.
In particular, our result applies to all forests.
As well, we give formulas for $\beta_{i,i+2}(I(G))$ for $i=0,\ldots,3$
for any simple graph $G$.
The formulas for $i=0,1$ were known, but our proofs are new.
We also give upper and lower bounds
for $\beta_{i,i+2}(I(G))$ for any simple graph based upon the numerical
information associated to the graph $G$.  
Using growth bounds on the Betti numbers of
lex ideals, we also obtain a purely graph theoretical 
result by providing 
a lower bound on the number of triangles in $G$.


\section{The graded Betti numbers of an edge ideal}

The main result of this section is to express $\beta_{i,j}(I(G))$
in terms of the dimensions of
the reduced simplicial homology groups of a simplicial complex
constructed from the graph $G$.

\subsection{Graphs and simplicial complexes}
Let $G$ denote a graph with vertex set $V_G =
  \{x_1,\ldots,x_n\}$ and edge set $E_G$.  
We shall say $G$ is {\it simple}
if $G$ has no loops or multiple
edges. A simple graph need not be connected.
The {\it degree} of a vertex
$x_i \in V_G$, denoted $\deg x_i$, is the number of edges incident to $x_i$.

If $S \subseteq V_G$, then the {\it induced subgraph} of $G$ on the
vertex set $S$, denoted $G_{S}$, is the
subgraph of $G$ such that every edge of $G$ contained in $S$
is an edge of $G_{S}$.  The {\it complement} of a graph $G$, denoted
by $\comp{G}$, is the graph whose vertex set is the 
same as $G$, but whose edge set is defined by
the rule: $\{x_i,x_j\} \in E_{G^c}$ if and only
if $\{x_i,x_j\} \not\in E_G$.  We shall let $\numbc(G)$ denote
the number of connected components of $G$.  Furthermore,
we let $\iv(G)$ denote the number of isolated vertices of $G$.
It is then clear that $\numbc(G) \geq \iv(G)$.

The {\it cycle} on $n$ vertices, denoted $C_n$, is the
graph with  $V_G = \{x_1,\ldots,x_n\}$ and $E_G
= \{\{x_1,x_2\},\{x_2,x_3\},\ldots,\{x_n,x_1\}\}$.   
A {\it forest} is any graph with no induced cycles; a {\it tree} is 
a connected forest.
The
{\it wheel} $W_n$ is the graph 
obtained by adding a vertex $z$ to $C_n$ and then
adjoining an edge between $z$ and every vertex $x_i \in V_{C_n}$.
Note that $W_n$ has $n+1$ vertices.
The {\it complete graph} on $n$ vertices, denoted $K_n$, 
is the graph with
the property that for $x_i \neq x_j \in V_{K_n}$, 
the edge $\{x_i,x_j\}\in E_{K_n}$.  The {\it complete bipartite graph},
denoted $K_{n,m}$, is the graph with vertex set $V_G = \{x_1,\ldots,x_n,y_1,
\ldots,y_m\}$ and edge set $E_G = \{\{x_i,y_j\} ~|~ 1 \leq i \leq n,~~
1 \leq j \leq m\}$.  We shall write $c_n(G),w_n(G),k_n(G),$ and $k_{n,m}(G)$
for the number of induced subgraphs of $G$ isomorphic to
$C_n,W_n,K_n$, and $K_{n,m}$, respectively.
 
A {\it simplicial complex} $\Delta$ on a set of vertices
$V = \{x_1,\ldots,x_n\}$ is a collection of subsets
of $V$ such that $\{x_i\} \in \Delta$ for $i = 1,\ldots,n$,
and for each $F \in \Delta$, if $G \subseteq F$, then
$G \in \Delta$.  Note that $\emptyset \in \Delta$. 
Associated to $\Delta$ is a square-free monomial
ideal $I_{\Delta}$ in $R = k[x_1,\ldots,x_n]$.  More precisely,
the ideal
\[I_{\Delta}  = \langle\{x_{i_1}x_{i_2}\cdots x_{i_t} ~|~ 
\{x_{i_1},x_{i_2},\ldots,x_{i_t}\}\not\in \Delta\}\rangle \]
is the {\it Stanley-Reisner ideal} of $\Delta$, and $I_{\Delta}$ is generated
by the square-free monomials that correspond to non-faces of $\Delta$.

If $G$ is any graph, the {\it clique complex} is the simplicial
complex $\Delta(G)$ where $F = \{x_{i_1},\ldots,x_{i_j}\} \in 
\Delta(G)$ if and only if $G_{F}$ is a complete graph.

If $G$ is a simple graph with edge ideal $I(G)$, then
the ideal $I(G)$ is generated by square-free monomials.
So $I(G)$ is also the defining ideal of a simplicial complex
$\Delta$, that is, $I(G) = I_{\Delta}.$  Specifically,
$I(G)$ is the Stanley-Reisner ideal associated to
the clique complex $\Delta = \Delta(\comp{G})$.  

\subsection{Graded Betti numbers of edge ideals}
Let $I$ be a homogeneous ideal of $R$.  The $ij$th
{\it graded Betti number} of $I$, denoted $\beta_{i,j}(I)$, 
is given by
\[\beta_{i,j}(I) = \dim_k \tor^R_i(I,k)_j.\]
If $I$ is also a monomial ideal, then $I$ also has an $\N^n$-graded
minimal free resolution. 
As a consequence $\tor^R_i(I,k)$ is also $\N^n$-graded,
and we have
\[\beta_{i,\alpha}(I) = \dim_k \tor^R_i(I,k)_{\alpha} 
~~\mbox{with $\alpha \in \N^n$}.
\]

When $I = I_{\Delta}$ is the Stanley-Reisner ideal
of a simplicial complex, Hochster \cite{Ho} 
described  the $\N^n$-graded Betti numbers of $I_{\Delta}$ in
terms of the dimensions of the reduced homology of the simplicial
complex $\Delta$.  To state the result, we require the following notation:
if 
$m$ is a  monomial of $R = k[x_1,\ldots,x_n]$ of multidegree $\alpha
\in \N^n$,
then we define
\[\tor_i^R(I_{\Delta},k)_m := \tor_i^R(I_{\Delta},k)_{\alpha}\]

\begin{proposition}[Hochster's Formula]\label{prop: hochster} 
Let $\Delta$ be a simplicial complex on vertex set
$\{x_1,\ldots,x_n\}$ and let $m$ be a monomial of $R$.
If $m$ is square-free, then
\[\dim_k \tor_i^R(I_{\Delta},k)_m =
\dim_k \widetilde{H}_{\deg(m)-i-2}(|m|,k)\]
where $\widetilde{H}_j(|m|,k)$ denotes the $j$th reduced homology
of the full subcomplex $|m|$ of $\Delta$ whose vertices correspond
to the variables dividing $m$.  If $m$ is not square-free,
then $\tor_i^R(I_{\Delta},k)_m$ vanishes.
\end{proposition}

By specializing this result to edge ideals, we obtain:

\begin{proposition}\label{prop: edgeidealbetti}
Let $G$ be a simple graph with edge ideal $I(G)$.  Then
\[\beta_{i,j}(I(G)) = \sum_{
\footnotesize
\begin{array}{c} 
S \subseteq V_G, ~~
|S| = j
\end{array}}
 \dim_k \widetilde{H}_{j-i-2}(\Delta(\comp{G_{S}}),k)
~~\mbox{for all $i,j\geq 0$.}\]
\end{proposition}

\begin{proof} We sketch out the main idea of the proof.
Let $\Delta = \Delta(G^c)$ be the simplicial complex defined by $I(G)$.
It follows from Proposition \ref{prop: hochster} that
\[\beta_{i,j}(I(G)) = \beta_{i,j}(I_{\Delta})
= \sum_{\footnotesize
\begin{array}{c} 
m \in M_j, ~~m ~\mbox{is square-free}
\end{array}}
\dim_k \widetilde{H}_{j-i-2}(|m|,k)\]
where $M_j$ consists of all the monomials of degree $j$ in $R$.
Since $\deg m = j$ and $m$ is square-free,
the variables that divide $m$ give a subset
$S \subseteq V_G$ of size $j$.  
Let $G_{S}$ denote the induced subgraph of $G$ on this vertex
set $S$, and let $\comp{G_{S}}$ denote its complement.  To finish the proof, 
it is enough to note that the full subcomplex $|m|$ of $\Delta(G^c)$ and 
$\Delta(\comp{G_{S}})$ are the same simplicial complex.
\end{proof}

\begin{remark}
A similar formula for $\beta_{i,j}(I(G))$ can be found in
\cite{K}.
While Proposition \ref{prop: edgeidealbetti} gives a formula
for all the numbers $\beta_{i,j}(I(G))$, the formula appears
difficult to apply since one has to compute 
the dimensions of all the homology groups 
$\widetilde{H}_{j-i-2}(\Delta(\comp{G_{S}}),k)$ as
$S$ varies over all subsets of $V_G$ of size $j$.
Moreover,
it is not clear how this formula relates to data associated to
$G$.
\end{remark}


\section{The Betti numbers in the linear strand of an edge ideal}
 
In this section we shall give some exact formulas
for the graded Betti numbers in
the linear strand of $I(G)$ when $G$ is a simple graph.

\begin{proposition}\label{prop: linearstrand}
Let $G$ be a simple graph with edge ideal $I(G)$.
Then 
\[
\beta_{i,i+2}(I(G)) = 
\sum_{\footnotesize
\begin{array}{c} 
S \subseteq V_G, ~~|S| = i+2
\end{array}}
 \left(\numbc(\comp{G_{S}})-1\right)
~~\mbox{for all $i\geq 0$.}\]
\end{proposition}

\begin{proof}
By Proposition \ref{prop: edgeidealbetti}
\[
\beta_{i,i+2}(I(G)) = \sum_{
\footnotesize
\begin{array}{c} 
S \subseteq V_G, ~~
|S| = i+2
\end{array}}
 \dim_k \widetilde{H}_{0}(\Delta(\comp{G_{S}}),k).
\]
Now $\dim_k \widetilde{H}_{0} (\Delta(\comp{G_{S}}),k) = 
\numbc(\Delta(\comp{G_{S}})) - 1$.
It now suffices to note that  
the simplicial complex $\Delta(\comp{G_{S}})$
and ${G^c_{S}}$ have
the same number of components.
\end{proof}

\begin{remark}
From this reformulation it is clear the numbers $\beta_{i,i+2}(I(G))$
do not depend upon the
characteristic of the field $k$.  See Katzman \cite{K} for
further results on the dependence of $\beta_{i,j}(I(G))$
upon the characteristic of $k$.
\end{remark}

\begin{remark}It will sometimes be convenient to write $\beta_{i,i+2}(I(G))$
as
{\footnotesize
\begin{equation}\label{splitequation}
\beta_{i,i+2}(I(G)) = 
\sum_{\footnotesize
\begin{array}{c} 
S \subseteq V_G,\\
|S| = i+2, \\
G^c_{S}\, \mbox{contains an}\\
\mbox{isolated vertex}\\
\end{array}}
 (\numbc(\compGS)-1) 
+
\sum_{\footnotesize
\begin{array}{c} 
S \subseteq V_G, \\
|S| = i+2,\\
G^c_{S}\, \mbox{contains no}\\
\mbox{isolated vertices}\\
\end{array}}
 (\numbc(\compGS)-1). 
\end{equation}}
\end{remark}
Although the formula in Proposition \ref{prop: linearstrand} 
is defined strictly in terms of the graph $G$,  
the formula has the disadvantage that one must sum over
all $S \subseteq V_G$
of size $i+2$.   A further reduction can be made if we impose
an extra hypothesis on $G$.  We say that $G$ has no
induced $4$-cycle if for every $S \subseteq V_G$ with $|S| =4$,
then $G_{S} \not\cong C_4$.

\begin{proposition}\label{prop: linearstrand_no_c4}
Let $G$ be a simple graph with edge ideal $I(G)$.  If $G$
has no induced $4$-cycles, then
\[\beta_{i,i+2}(I(G)) = \sum_{v \in V_G} \binom{\deg v}{i+1}
- k_{i+2}(G) 
~~\mbox{for all $i \geq 0$.}\]
Furthermore, the above formula holds for all simple graphs $G$ if $i=0$ or $1$.
\end{proposition}

\begin{proof} 
The hypothesis that $G$ contains no induced $4$-cycles is the same
as saying that for any subset $S$ of the vertices, at most
one connected component of $G_{S}^c$ is larger than a vertex.
That is, if $\compGS$ has two 
connected components, each containing an edge, then letting 
$S'$ be the set of the four vertices in these two edges we see that 
$G_{S'}$ is an induced $4$-cycle of $G_{S}$.
Therefore, to count $\numbc(\compGS)-1$ we can simply count the number of 
isolated vertices in $\compGS$, the
``$-1$'' term being taken care of by the component which isn't a vertex;  
we thus overcount by one whenever 
$\compGS$ consists completely of isolated vertices.

For any vertex $v$, the number of subsets $S$ of size $i+2$ containing 
$v$ such that $v$ is an isolated vertex in $\compGS$ is 
$\binom{\deg v}{i+1}$.  
To take care of the overcount, note that $\compGS$ consists of isolated
vertices
exactly when $G_{S}$ is a complete graph on $i+2$ vertices.  Subtracting the 
number of times this happens gives the formula above. 

Now let $G$ be any simple graph. To compute $\beta_{i,i+2}(I(G))$
when $i = 0$ or $1$, we need to count the number of 
connected components of $\comp{G_{S}}$ when $|S| =2$ or $3$.
But for any simple graph on two or three vertices, at most
one connected component can be larger than a vertex.  The proof
is now the same as the one given above.
\end{proof}

\begin{remark}  It follows from the construction of $I(G)$ 
that $\beta_{0,2}(I(G)) = \#E_G$.  Coupling
this result with the above formula gives the identity
\[\#E_G = \sum_{v \in V_G} \deg v  - k_2(G).\]
But $k_2(G) = \#E_G$ since each edge is isomorphic to $K_2$.
We thus recover the well known Degree-Sum Formula which 
states that $2\#E_G = \sum_{v \in V_G} \deg v$.
When $i=1$, we have the formula $\beta_{1,3}(I(G)) = \sum_{v \in V_G} 
\binom{\deg v}{2} - k_3(G)$.  Since $k_3(G)$ is simply the number
of triangles in $G$, we recover the formula of \cite{EV}.
\end{remark}

Since forests do not have cycles, we get an immediate corollary.

\begin{corollary} Let $G$ be a forest with edge ideal $I(G)$.  Then
$\beta_{0,2}(I(G)) = \#E_G$, and 
\[\beta_{i,i+2}(I(G)) = \sum_{v \in V_G} \binom{\deg v}{i+1}
~~\mbox{for all $i \geq 1$.}\]
\end{corollary}

\begin{remark}
Zheng \cite{Z} studied the graded Betti numbers of a special
class of simplicial complexes that generalized the notion
of a tree (due to Faridi \cite{F}).  
Proposition 3.3 of \cite{Z} gives an explicit formula for
the graded Betti numbers in the linear strand of a $d$-dimensional
pure tree $\Delta$ connected in codimension 1.  When $d = 1$,
$\Delta$ is simply a connected tree.  Zheng's formula agrees
with ours in this case.
\end{remark}

We now provide formulas for $\beta_{2,4}(I(G))$ and $\beta_{3,4}(I(G))$
for any simple graph $G$.  The formula for $\beta_{2,4}(I(G))$ verifies the
conjecture of \cite{EV}.

\begin{proposition}\label{exact form}
Let $G$ be any simple graph with edge ideal $I(G)$.  Then
\begin{eqnarray*}
\beta_{2,4}(I(G)) &=& \sum_{v\in V_G} \binom{\deg v }{3} - k_4(G) + 
k_{2,2}(G)~~\mbox{and} \\
\beta_{3,5}(I(G)) & = &\sum_{v \in V_G}\binom{\deg v}{4} - k_5(G) + k_{2,3}(G)
+ w_4(G) + d(G)
\end{eqnarray*}
where  $d(G)$ 
is the number of induced 
subgraphs of $G$ isomorphic to the graph $D$ given below:

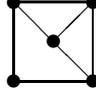
\begin{picture}(-100,40)(-155,0)
\put(0,0){\circle*{5}}
\put(0,0){\line(1,0){30}}
\put(0,0){\line(0,1){30}}
\put(0,30){\circle*{5}}
\put(0,30){\line(1,0){30}}
\put(30,0){\circle*{5}}
\put(30,0){\line(0,1){30}}
\put(30,30){\circle*{5}}
\put(15,15){\circle*{5}}
\put(0,30){\line(1,-1){30}}
\put(15,15){\line(1,1){15}}
\end{picture}
\end{proposition}

\begin{proof}
Suppose $i=2$.  We will compute both sums in $(\ref{splitequation})$
to find the value of $\beta_{2,4}(I(G))$.  Suppose $S \subseteq V_G$ is 
such that $G_{S}^c$ has an isolated vertex.  
The remaining three vertices of $G_{S}^c$ must form a simple graph.
By considering all possible simple graphs on three vertices,
we will have
\[\numbc(G_{S}^c) - 1 = \iv(G_{S}^c) \]
except when $G_{S}^c$ is the graph of 4 isolated vertices,
that is, $G_{S} \cong K_4$, in which case we have
\[\numbc(G_{S}^c) -1  = \iv(G_{S}^c) -1.\]
The first sum in the expression for $\beta_{2,4}(I(G))$ therefore equals
{\small
\begin{eqnarray*}
\sum_{\footnotesize
\begin{array}{c} 
S \subseteq V_G,~~|S| = 4, \\
G^c_{S}\, \mbox{contains an isolated vertex}\\
\end{array}} 
\iv(G_{S}^c)  - \sum_{G_{S} \cong K_4} 1  
 =  \sum_{v \in V_G} \binom{\deg v}{3} - k_4(G).
\end{eqnarray*}}

On the other hand, suppose $S \subseteq V_G$ is such that
$G_{S}^c$ has no isolated vertex.  We can further assume that
$G_{S}^c$ is disconnected, because if $G_{S}^c$
is connected, no contribution is made to $\beta_{2,4}(I(G))$ for
this subset.  But then the only possibility for $G_{S}^c$
is two disjoint edges, or equivalently, $G_{S} \cong K_{2,2}$.  So the
second sum in $(\ref{splitequation})$ must equal $k_{2,2}(G)$,
thus giving the desired formula for $\beta_{2,4}(I(G))$.

Suppose now that $i=3$.  If $S \subseteq V_G$ is such that $|S| = 5$
and $G_{S}^c$ has an isolated vertex, then the remaining
four vertices form a simple graph.  By considering all possible simple
graphs on four vertices, we have
\[
\numbc(G_{S}^c) - 1 =
\left\{
\begin{array}{ll}
\iv(G_{S}^c) + 1 & \mbox{if $G_{S} \cong W_4$} \\
\iv(G_{S}^c) - 1 & \mbox{if $G_{S} \cong K_5$} \\
\iv(G_{S}^c) & \mbox{otherwise}
\end{array}\right..\]
Thus
\[
\sum_{\footnotesize
\begin{array}{c} 
S \subseteq V_G,~~|S| = 5, \\
G^c_{S}\, \mbox{contains an isolated vertex}\\
\end{array}}
 (\numbc(\compGS)-1) = \sum_{v \in V_G} \binom{\deg v}{4} - k_5(G) + w_4(G).\]

Now suppose $S$ is such that $|S| = 5$ but $G_{S}^c$ has
no isolated vertices.  Since such an $S$ will only contribute
to $\beta_{3,5}(I(G))$ if $G_{S}^c$ is disconnected, we
can further assume that $S$ is such that $G_{S}^c$ is disconnected.
By considering all simple graphs on 5 vertices, we find
there are only two possibilities for $S$:  either $G_{S} \cong K_{2,3}$, 
or $G_{S}$ is isomorphic to
the graph $D$.  In both cases, $\numbc(G_{S}^c) - 1 = 1$.  Hence
\[
\sum_{\footnotesize
\begin{array}{c} 
S \subseteq V_G, ~~|S| = 5,\\
G^c_{S}\, \mbox{contains no isolated vertices}\\
\end{array}}
 (\numbc(\compGS)-1) = k_{2,3}(G) + d(G). \]
If we substitute these values into $(\ref{splitequation})$, we get the desired
formula for $\beta_{3,5}(I(G))$.
\end{proof}

For any  homogeneous ideal $I \subseteq R$ generated in degree 
$d$,  we say $I$ has a {\it linear resolution} if $\beta_{i,j}(I) = 0$
for any $j \neq i + d$.  In other words, the only non-zero graded
Betti numbers of $I$ are those in the linear strand.  When
$I = I(G)$ for some simple graph $G$, Fr\"oberg
\cite{Fr} gave a characterization
of those ideals having a linear resolution in terms of chordal
graphs.  
We say $G$ is a {\it chordal graph} if every cycle of length
$n > 3$ has a chord.  Here, if $\{x_1,x_2\},\ldots,\{x_n,x_1\}$
are the $n$ edges of a cycle of length $n$, we say the
cycle has a chord in $G$ if there exists vertices $x_i,x_j$ 
such that $\{x_i,x_j\}$ is also an edge of $G$, 
but $j \not\equiv i \pm 1 \pmod{n}$.

\begin{proposition}\label{linearresolution}
Let $G$ be a graph with edge ideal $I(G)$.  The edge ideal $I(G)$
has a linear resolution if and only if $\comp{G}$ is a chordal graph.
\end{proposition}

Using this criterion and the results of this section we can completely 
describe 
the graded minimal free resolutions of $I(G)$ when $G = K_n$ or $K_{n,m}$.
These results were proved via different means in \cite{J}.

\begin{example}(The resolution of $I(K_n)$)  The complement
of $K_n$ is simply $n$ isolated vertices, and hence,
$I(K_n)$ must have a linear resolution by Proposition \ref{linearresolution}.
Moreover, no induced subgraph of $K_n$
will be a $4$-cycle.  Proposition \ref{prop: linearstrand_no_c4}
then gives
$\beta_{i,i+2}(I(K_n)) = \sum_{v \in V_{K_n}} 
\binom{\deg v}{i+1} - k_{i+2}(K_n) ~\mbox{for all $i\geq 0$}.$
Each vertex of $K_n$ has degree $n-1$, and $K_n$ has $\binom{n}{i+2}$
subgraphs isomorphic to $K_{i+2}$.  The above expression
therefore reduces to
\[\beta_{i,i+2}(I(K_n)) = (i+1)\binom{n}{i+2}~\mbox{for all $i\geq 0$}.\]
\end{example}

\begin{example}(The resolution of $I(K_{a,b})$)
Let $G = K_{a,b}$ be a complete bipartite graph.  We write
the vertex set of $G$ as $V_G = \{x_1,\ldots,x_a,y_1,\ldots,y_b\}$
so that $E_G = \{\{x_i,y_j\} ~|~ 1 \leq i \leq a, ~1\leq j \leq b\}$.
For all $a,b \geq 1$, the complement of $G$ is the 
disjoint union of $K_a$ and $K_b$.  Since this graph
has no induced cycle of length $\geq 4$, the resolution 
of $I(G)$ is linear by Proposition \ref{linearresolution}.

Because $\comp{G} = K_a \cup K_b$,
for any $S \subseteq V_G$ with $|S| = i+2$, we have
\[
\numbc(\comp{G_{S}}) = 
\left\{
\begin{array}{ll}
2 & \text{if}~S \cap\{x_1,\ldots,x_a\} \neq \emptyset
~\text{and}~S \cap\{y_1,\ldots,y_b\} \neq \emptyset \\
1 & \mbox{otherwise}
\end{array}
\right..
\]
By Proposition \ref{prop: linearstrand}, to determine
$\beta_{i,i+2}(I(G))$ it therefore suffices to 
count the number of subsets $S \subseteq V_G$ with
$|S| = i+2$ and $\numbc(\comp{G_{S}}) =2$.

There are $\binom{a+b}{i+2}$ subsets of $V_G$ that contain
$i+2$ distinct vertices.  Furthermore, $\binom{a}{i+2}$
of these subsets must contain only vertices among
$\{x_1,\ldots,x_a\}$; similarly, $\binom{b}{i+2}$
of these subsets contain only vertices among $\{y_1,\ldots,y_b\}$.
It thus follows that
\[\beta_{i,i+2}(I(K_{a,b})) = \binom{a+b}{i+2} -
\binom{a}{i+2} - \binom{b}{i+2} ~~\text{for all $i \geq 0$}\]
since the expression on the right hand side counts the number
of subsets $S \subseteq V_G$ with $|S| = i+2$ and $S$ contains
at least one $x_i$ and one $y_j$ vertex.
\end{example}

\section{Bounds on Betti numbers in the linear strand}

In
this section we provide some bounds on $\beta_{i,i+2}(I(G))$ for
any graph $G$.

\begin{proposition}\label{prop: linearstrand_lower_bound}
Let $G$ be a simple graph with edge ideal $I(G)$.
Then for all $i \geq 0$
{\small
\[\beta_{i,i+2}(I(G)) \geq \sum_{v \in V_G} \binom{\deg v}{i+1}
- k_{i+2}(G) + k_{2,i}(G) + k_{3,i-1}(G) + 
\cdots + k_{\lfloor \frac{i+2}{2}\rfloor,\lceil \frac{i+2}{2}\rceil}(G).\]}
\end{proposition}

\begin{remark}
When $G$ has no induced $4$-cycle we have  $k_{j,i+2-j}(G) = 0$ for
$j = 2,\ldots,\lfloor \frac{i+2}{2}\rfloor$.  In this case 
$\beta_{i,i+2}(I(G))$ agrees with the lower bound
by Proposition \ref{prop: linearstrand_no_c4}.
The numbers $\beta_{i,i+2}(I(G))$ are also given by lower bound
when $i = 0,1,$ or $2$.
\end{remark}

\begin{proof} 
We will use the formula for $\beta_{i,i+2}(I(G))$ as given
in equation $(\ref{splitequation})$.
Suppose that $S \subseteq V_G$ is such
that $G_{S}^c$ contains an isolated vertex.  If $G_{S}^c$
consists entirely of isolated vertices, i.e.,
$G_{S} \cong K_{i+2}$, then 
\[\numbc(G_{S}^c) - 1 = \iv(G_{S}^c) -1. \]
On the hand, if $G_{S}^c$ has at least one component consisting of an
 edge, then
\[\numbc(G_{S}^c) - 1 \geq  \iv(G_{S}^c) \]
since the ``-1'' takes care of the component containing the edge.
It thus follows that
\[\sum_{\footnotesize\begin{array}{c} 
S \subseteq V_G,\\
|S| = i+2, \\
G^c_{S}\, \mbox{contains an}\\
\mbox{isolated vertex}\\
\end{array}}
 (\numbc(\compGS)-1) \geq \sum_{\footnotesize
\begin{array}{c} 
S \subseteq V_G,\\
|S| = i+2, \\
G^c_{S}\, \mbox{contains an}\\
\mbox{isolated vertex}\\
\end{array}}
 \iv(G_{S}^c) - \sum_{
\footnotesize \begin{array}{c} 
S \subseteq V_G \\
G_{S} \cong K_{i+2}
\end{array}} 1.\]
For any vertex $v$, the number of subsets $S$ of size $i+2$ containing $v$
such that $v$ is an isolated vertex in $G_{S}^c$ is 
$\binom{\deg v}{i+1}$.  We thus obtain the bound
\begin{equation}\label{bound1}
\sum_{\footnotesize
\begin{array}{c} 
S \subseteq V_G,~~
|S| = i+2, \\
G^c_{S}\, \mbox{contains an isolated vertex}\\
\end{array}}
 (\numbc(\compGS)-1) \geq 
\sum_{v \in V_G} \binom{\deg v}{i+1} - k_{i+2}(G).
\end{equation}

Suppose that  $S \subseteq V_G$ is such that
$G_{S}^c$ contains no isolated vertices.  If $G_{S}^c$ is connected,
then this subset $S$ does not contribute to the value of $\beta_{i,i+2}(I(G))$
because $\numbc(G_{S}^c) - 1 = 0$.
So, we can assume that $S$ is such that $G_{S}^c$ has at 
least two connected components.  

Now if $S$ is such that
$|S| = i+2$ and 
$G_{S} \cong K_{a,i+2-a}$ for some 
$a = 2,\ldots,\lfloor \frac{i+2}{2} \rfloor$, then $G_{S}^c$ has exactly
two components, and no component is an isolated vertex.  
Thus, each such $S$ contributes exactly 1 to the value of 
$\beta_{i,i+2}(I(G))$, and 
hence
\begin{equation}\label{bound2}
\sum_{\footnotesize
\begin{array}{c} 
S \subseteq V_G ~~
|S| = i+2,\\
G^c_{S}\, \mbox{contains no isolated vertices}\\
\end{array}}
 (\numbc(\compGS)-1) \geq k_{2,i}(G) + \cdots 
+ k_{\lfloor\frac{i+2}{2}\rfloor,
\lceil\frac{i+2}{2}\rceil}(G).
\end{equation}
By now combining the bounds (\ref{bound1}) and (\ref{bound2}),
we get the desired conclusion.
\end{proof}

Recall that the elements of $\N^n$ can be given a total ordering
using the lexicographical order defined by $(a_1,\ldots,a_n)
> (b_1,\ldots,b_n)$ if $a_1 = b_1, \ldots, a_{i-1} = b_{i-1}$
but $a_i > b_i$.  This induces an ordering on the monomials
of $R$:  $x_1^{a_1}\cdots x_{n}^{a_n} >_{lex} x_1^{b_1}\cdots x_n^{b_n}$
if $(a_1,\ldots,a_n) > (b_1,\ldots,b_n)$.  A monomial
ideal $I$ is a {\it lex ideal} if for each $d \in \N$,
a basis for $I_d$ is the $\dim_k I_d$ largest monomials of degree $d$
with respect to the lexicographical ordering. 

\begin{proposition}\label{prop: upperbound}
Let $G$ be a simple graph and let
$\{m_1,\ldots,m_{\#E_G}\}$ be the $\#E_G$ largest monomials
of degree 2 in $R$ with respect to the lexicographical 
ordering.  Then
\[\beta_{i,i+2}(I(G)) \leq \sum_{t=1}^{\#E_G} \binom {u_t -1}{i}\]
where $u_t$ is the largest index of a variable dividing $m_t$. 
\end{proposition}
 
\begin{proof}
Let $L$ be the lex ideal 
with the property that  $H_{R/L}$, the Hilbert function of $R/L$, 
is the same as the Hilbert function of  $R/I(G)$.  By 
the Bigatti-Hulett-Pardue Theorem (see \cite{B,H,P}),  we have 
\[\beta_{i,j}(I(G)) \leq \beta_{i,j}(L) ~~\mbox{for all $i,j \geq 0$}.\]
In particular, $\beta_{i,i+2}(I(G)) \leq \beta_{i,i+2}(L)$ for
all $i\geq 0$.  Since $L$ is generated in degrees greater than or
equal to $2$, it follows that
$\beta_{i,i+2}(L) = \beta_{i,i+2}(L_2) ~~\mbox{for all $i \geq 0$}$
where $L_2$ is the ideal generated by  all the degree 2 monomials
in $L$.  Moreover, since $H_{R/I(G)}(2) = H_{R/L}(2)$, $L$ has 
$\#E_G$ monomials of degree 2.  Thus $L_2 = \langle m_1,\ldots,m_{\#E_G}
\rangle$
where the $m_i$s are the $\#E_G$ largest monomials of degree 2
with respect to the lexicographical ordering.  Note that
$L_2$ is also a lex ideal.  

To get the desired conclusion, we use the Eliahou-Kervaire resolution 
(see \cite{EK})
of Borel-fixed ideals
(a class of ideals which includes the
lexicographical ideals) to compute $\beta_{i,i+2}(L_2)$:
\[\beta_{i,i+2}(L_2) =   \sum_{t=1}^{\#E_G} \binom {u_t -1}{i}
~~\mbox{for all $i \geq 0$}\]
where $u_t$ is the largest index of a variable dividing $m_t$. 
\end{proof}

\begin{remark}
We can give an explicit description of the $\#E_G$ largest
monomials of degree 2.  Let $j$ be
the integer such that
\[n+ (n-1) + \cdots + n-(j-1) < \#E_G \leq n + (n-1) +\cdots + n-j,\]
and let $l = \#E_G - \sum_{i=0}^{j-1} (n-i).$
(If $j =0$, then we take $\sum_{i=0}^{j-1} (n-i) = 0$.)
Then the $\#E_G$ largest monomials of degree 2 are
\[\{x_1^2, x_1x_2,\ldots,x_1x_n,x_2^2,x_2x_3,\ldots,x_2x_n,x_3^2,\ldots, x_{j+1}^2,\ldots,
x_{j+1}x_{j+l}\}.\]
Notice that since we are only interested in the indices,
which only depend upon $\#E_G$ and $\#V_G$, 
Proposition \ref{prop: upperbound} can be written 
without reference to monomials.
\end{remark}

The bound in the previous result is not related to the
numerical information of the graph $G$, but is based upon results
describing the growth of Betti numbers of monomial ideals.  By playing
this result against the exact formulas of the previous section, 
we can derive bounds on information describing the graph.  For example,
as a corollary, we derive a lower bound for the number
of triangles in a graph.  

\begin{corollary} Let $G$ be a simple graph, and let $j$ and $l$
be defined from $\#E_G$ as above.  Then
\[
\max\left\{0, \sum_{v \in V_G} \binom{\deg v}{2}  -j\binom{\#V_G}{2} +
\binom{j}{j-3} - (j+\cdots + (j+l-1))\right\}
\leq k_3(G).\]
\end{corollary}
\begin{proof}
Let  $\{x_1^2, x_1x_2,\ldots, x_{j+1}^2,\ldots,
x_{j+1}x_{j+l}\}$
be the $\#E_G$ largest monomials of degree 2 with respect to lexicographical
ordering.
By combining Propositions \ref{prop: linearstrand_no_c4} and
\ref{prop: upperbound} we get
\[\beta_{1,3}(I(G)) = \sum_{v \in V_G} \binom{\deg v}{2} - k_3(G) \leq
\sum_{t=1}^{\#E_G} u_t - 1 \]
where $u_t$ is the largest index of a variable dividing $m_i$.
From our description of the monomials, it is now a calculation
to check that 
\[\sum_{t=1}^{\#E_G} u_t - 1 = j\binom{\#V_G}{2} -
\binom{j}{j-3} + (j+\cdots + (j+l-1)).\]
(We have used the fact that $n = \#V_G$.)
After substituting this result into our inequality for
$\beta_{1,3}(I(G))$ and rearranging, we the get desired bound.
\end{proof}

\begin{remark}
The above result is far from being the best lower bound
for the number of triangles in a graph.  However, the above results
nicely illustrate that a better 
understanding of the behaviour of Betti numbers of monomial ideals will give
new algebraic tools to prove graph theoretical results.
\end{remark}

\noindent
{\it Acknowledgments.}  We would like to thank Huy T\`ai H\`a and 
David Gregory for their comments and some fruitful discussions.
Both authors would like to acknowledge the financial support of NSERC while
working on this project.  The second author
would also like to thank Queen's University
and the Universit\`a di Genova for their hospitality while working
on this project.


\end{document}